\documentclass[12pt]{amsart}

\usepackage{amsmath}
\usepackage{latexsym}
\usepackage{amssymb}
\usepackage{amsfonts}
\usepackage{color}

\usepackage{graphics}

\renewcommand{\proof}{\par\noindent{\it Proof.\ \ }}
\def\qed{\ifmmode\square\else\nolinebreak\hfill
$\square$\fi\par\vskip12pt}

\renewcommand{\proof}{\par\noindent{\it Proof.\ \ }}
\def\qed{\ifmmode\square\else\nolinebreak\hfill
$\square$\fi\par\vskip12pt}

\def\l{\langle} \def\r{\rangle} 
\def\div{\,\big|\,} 

 \def\ZZ{\mathbb Z}

 \def\Cay{{\sf Cay}} \def\Cos{{\sf Cos}}

\def\D{{\rm D}} \def\S{{\rm S}} 
 \def\M{{\sf M}}

 \def\Aut{{\sf Aut}}

\def\ASL{{\sf ASL}}

  \def\K{{\sf K}}

\def\Mult{{\sf Mult}}

\def\Ga{{\it \Gamma}}

\def\Del{{\it\Delta}}
\def\Ome{{\it\Omega}}

\def\a{\alpha}

 \def\PSp{{\rm PSp}}

\def\A{{\rm A}}

\def\PSL{{\rm PSL}}
\def\GL{{\rm GL}} \def\SL{{\rm SL}}

 \def\PSU{{\rm PSU}}

 \def\F{{\rm F}} \def\D{{\rm D}}

\newtheorem{theorem}{Theorem}[section]
\newtheorem{lemma}[theorem]{Lemma}%
\newtheorem{proposition}[theorem]{Proposition}%
\newtheorem{example}[theorem]{Example}%

    \def\ZZ{\mathbb Z}

\begin{document}

\title[Cayley graphs on non-abelian simple groups]{Arc-transitive Cayley graphs
on non-ableian simple groups with soluble vertex stabilizers and valency seven}
\thanks{1991 MR Subject Classification 20B15, 20B30, 05C25.}
\thanks{This paper was partially supported by the National Natural Science Foundation of China (11461004,11231008).}

\author[Pan]{Jiangmin Pan}
\address{J. M. Pan\\
School of Statistics and Mathematics\\
Yunnan University of Finance and Economics\\
Kunming \\
P. R. China}
\email{jmpan@ynu.edu.cn}

\author[Yin]{Fugang Yin}
\address{F. G. Yin\\
School of Statistics and Mathematics\\
Yunnan University of Finance and Economics\\
Kunming \\
P. R. China}

\author[Ling]{Bo Ling}
\address{B. Ling, School of Mathematics and Computer Science\\
Yunnan Minzu University\\
Kunming, Yunnan, P. R. China}
\email{bolinggxu@163.com}

\curraddr{}


\maketitle

\begin{abstract}
In this paper, we study arc-transitive Cayley graphs on non-abelian simple groups with
soluble stabilizers and valency seven.
Let $\Ga$ be such a Cayley graph on a non-abelian simple group $T$. It is proved that
either $\Ga$ is a normal Cayley graph or
$\Ga$ is $S$-arc-transitive,
with $(S,T)=(\A_n,\A_{n-1})$
and $n=7,21,63$ or $84$;
and, for each of these four values of $n$, there really exists
arc-transitive $7$-valent non-normal
Cayley graphs on $\A_{n-1}$ and specific examples are constructed.

\end{abstract}

\qquad {\textsc k}{\scriptsize \textsc {eywords.}} {\footnotesize
arc-transitive graph, Cayley graph, non-abelian simple group, coset graph}

\section{Introduction}
Throughout the paper, all groups are finite,
and all graphs are finite, undirected and simple.

For a graph $\Ga$, denote by $V\Ga$ its vertex set,
and $\Aut\Ga$ its full automorphism
group. If there is an automorphism group $X\le\Aut\Ga$
which is transitive on the vertex set or the arc set
of $\Ga$,
then $\Ga$ is called {\it $X$-vertex-transitive} or
{\it $X$-arc-transitive}, respectively.
An arc-transitive graph is also called {\it symmetric}.
Let $s$ be a positive integer. An $s$-{\it arc} of $\Ga$
is a sequence $\a_0,\a_1,\dots,\a_s$ of $s+1$ vertices such that
$\a_{i-1},\a_i$ are adjacent for $1\le i\le s$ and $\a_{i-1}\not=\a_{i+1}$
for $1\le i\le s-1$. Then
$\Ga$ is called {\it $(X,s)$-arc-transitive} if $X\le\Aut\Ga$ is transitive on the set of
$s$-arcs of $\Ga$. If $X$ is transitive on the set of
$s$-arcs but not transitive on the set of $(s+1)$-arcs of $\Ga$,
then $\Ga$ is called {\it $(X,s)$-transitive},
and $\Ga$ is simply called {\it $s$-transitive} while $X=\Aut\Ga$.

Cayley graphs provide a rich source of transitive graphs.
A graph $\Ga$ is called a {\it Cayley graph}
on a group $G$ if there is a subset $S\subseteq G\setminus\{1\}$,
with $S=S^{-1}:=\{s^{-1}\mid s\in S\}$,
such that its vertex set equals $G$
and two vertices $g$ and $h$ are adjacent if and
only if $hg^{-1}\in S$. We denote this Cayley graph by $\Cay(G,S)$.
It is well known that the right regular
presentation
$$\hat G=\{\hat g\mid \hat g:~x\mapsto xg,\ \mbox{for all $g,x\in G$}\}$$
of $G$ is a subgroup of $\Aut\Ga$.
If $\hat G$ is normal in $\Aut\Ga$,
then $\Ga$ is called a {\it normal Cayley graph} on $G$.
Clearly, $\hat G$ is isomorphic to $G$,
for convenience, we will write $\hat G$ as $G$ in the following.

Due to the fundamental significance of the non-abelian simple groups in group theory,
studying Cayley graphs on non-abelian simple groups
is naturally an interesting topic in algebraic graph theory.
In 1996, by using a so called `dual action', Li \cite{Li96}
proved that a connected cubic arc-transitive Cayley graph on a non-abelian simple group
$T$ is either normal or $T=\A_5,\A_{11},\A_{23},\A_{47},\PSL(2,11),\M_{11}$
or $\M_{23}$.
Li's result was extended by Xu et al. \cite{Xu05,Xu07}
by proving that if $T\ne\A_{47}$,
then $\Ga$ is normal, and
there are exactly two non-isomorphic
cubic arc-transitive
non-normal Cayley graphs on $\A_{47}$.
In 2002, Fang, Praeger and Wang \cite{FPW02}
developed a theory for investigating the automorphism
groups of Cayley graphs on non-abelian simple groups,
which is then used to characterize tetravalent
and pentavalent arc-transitive Cayley graphs on non-abelian simple groups
by \cite{FLX04,FMW11}. Recently, the result in \cite{FMW11}
is extended by \cite{DFZ17} by characterizing
pentavalent arc-transitive graphs admitting
a vertex-transitive non-abelian simple group.

The main purpose of this paper is to present the following classification regarding
arc-transitive Cayley graphs on non-abelian simple groups
with valency seven and soluble stabilizers.

\begin{theorem}\label{Thm-1}
Let $\Ga$ be a connected arc-transitive Cayley graph on
a non-abelian simple group $T$
with valency seven and soluble stabilizer.
Then one of the following statements holds.

\begin{itemize}
\item[(1)] $\Ga$ is a normal Cayley graph.
\item[(2)] $\Ga$ is a non-normal Cayley graph, and either
\begin{itemize}
\item[(i)] $\Ga$ is
$(S,1)$-transitive
with $(S,T)=(\A_n,\A_{n-1})$ and $n=7,21$ or $63$,~or
\item[(ii)] $\Ga$ is $(\A_{84},2)$-transitive
with $T=\A_{83}$.
\end{itemize}
\end{itemize}

Further, for each $n\in\{7,21,63,84\}$, there
exists non-normal Cayley graph on $\A_{n-1}$ satisfying parts $(i)$ and $(ii)$.

\end{theorem}

{\noindent\bf Remarks.} (1) With certain nice properties of normal Cayley graphs,
graphs in part (1) can be well-characterized, see \cite{Li-ETCay}.

(2) Graphs in part (2) can be completed determined by using coset graphs and Magma
computational package \cite{Magma}.

(3) Specific examples which satisfy parts (i) are (ii) are constructed in Section 5.

\vskip0.1in
After this introductory section, we introduce some preliminary results
in Section 2. Then by proving some technical lemmas in Section 3,
we complete the proof of Theorem~\ref{Thm-1} in Section 4.
In the final Section 5, we present some specific examples of graphs
satisfying parts (i) and (ii) of Theorem~\ref{Thm-1}.

\section{Preliminaries}

In this section, we give some preliminary results
that will be used in the subsequent discussions.

We begin with some notational conventions used throughout this paper.
For a positive integer $n$, we use $\ZZ_n$,
$\D_{n}$ (if $n$ even) and $\F_n$
to denote the cyclic group, the dihedral group and Frobenius (but not dihedral) group of order $n$,
respectively.
For a group $G$ and a subgroup $H$ of $G$, denote by
$N_G(H)$ and $C_G(H)$
the normalizer and the centralizer of $H$ in $G$ respectively.
Given two groups $N$ and $H$, denote by $N\times H$ the direct product of $N$ and $H$,
by $N.H$ an extension of $N$ by $H$, and
if such an extension is split, then we write $N:H$ instead of $N.H$.

The following is a so called `N/C' theorem.

\begin{lemma}\label{N/C}{\rm(\cite[Ch. I, 1.4]{Huppert})}
Let $G$ be a group and $H$ a subgroup of $G$.
Then $N_G(H)/C_G(H)\lesssim\Aut(H)$.
\end{lemma}

The next is {\it Frattini argument} on transitive permutation groups.

\begin{lemma}\label{Frattini}{\rm(\cite[P. 9]{DM})}
Let $G$ be a transitive permutation group on $\Ome$,
$H$ a subgroup of $G$ and $\a\in\Ome$.
Then $H$ is transitive if and only if $G=HG_{\a}$.
\end{lemma}

The vertex stabilizers of $7$-valent arc-transitive graphs
were determined independently by \cite[Theorem 1.1]{GLH16}
and \cite[Theorem 3.4]{LLW16}.

\begin{lemma}\label{Stabilizer}
Let $\Ga$ be a connected $7$-valent $(X,s)$-transitive graph, where $X\le\Aut\Ga$ and $s\ge 1$.
Then $s\le3$ and one of the following statements holds, where $\a\in V\Ga$.

\begin{itemize}
\item[(1)] If $X_{\a}$ is soluble, then $|X_{\a}|\div 252$
and the couple $(s,X_{\a})$ is listed in the following table.

\[\begin{array}{c|c|c|c} \hline
s & 1 & 2 & 3 \\ \hline
X_{\a} &  \ZZ_7,~\D_{14},~\F_{21},~\D_{14}\times\ZZ_2,~\F_{21}\times\ZZ_3& \F_{42},
~\F_{42}\times\ZZ_2,~\F_{42}\times\ZZ_3& \F_{42}\times\ZZ_6  \\ \hline
\end{array}\]

\vskip0.1in
\item[(2)] If $X_{\a}$ is insoluble, then $|X_{\a}|\div 2^{24}\cdot 3^4\cdot 5^2\cdot 7$,
and the triple $(s,X_{\a},|X_{\a}|)$ is listed in the following table.

\[\begin{array}{c|c|c} \hline
s & 2 & 3  \\ \hline
X_{\a} & \PSL(3,2), \ASL(3,2),&  \PSL(3,2){\times}\S_4,\A_7{\times}\A_6,\S_7{\times}\S_6,(\A_7{\times}\A_6){:}\ZZ_2,   \\
       & \ASL(3,2){\times}\ZZ_2,\A_7,\S_7 &
       \ZZ_2^6{:}(\SL(2,2){\times}\SL(3,2)),[2^{20}]{:}(\SL(2,2){\times}\SL(3,2))  \\ \hline
       |X_{\a}|&2^3{\cdot}3{\cdot}7,2^6{\cdot}3{\cdot}7,&2^6{\cdot}3^2{\cdot}7,~
       2^6{\cdot}3^4{\cdot}5^2{\cdot}7,~2^8{\cdot}3^4{\cdot}5^2{\cdot}7,~2^7{\cdot}3^4{\cdot}5^2{\cdot}7,\\
         & 2^7{\cdot}3{\cdot}7,2^3{\cdot}3^2{\cdot}5{\cdot}7,
         2^4{\cdot}3^2{\cdot}5{\cdot}7& 2^{10}{\cdot}3^2{\cdot}7,~~2^{24}{\cdot}3^2{\cdot}7   \\ \hline
\end{array}\]
\end{itemize}
\end{lemma}

\vskip0.1in
A group $G$ is called {\it perfect} if $G = G'$, the commutator subgroup,
and an extension $G=N.H$ is called a {\it central extension} if $N\subseteq Z(G)$, the center of $G$.
If a group $G$ is perfect and $G/Z(G)$ is isomorphic to a simple group $T$, then $G$ is
called a {\it covering group} of $T$.
Schur \cite{Schur} showed that a simple (and, more generally,
perfect) group $T$ possesses a ¡®universal¡¯ covering group $G$ with the property that every
covering group of $T$ is a homomorphic image of $G$,
in this case, the center $Z(G)$ is called the {\it Schur multiplier} of $T$,
denoted by $\Mult(T)$.
The Schur multipliers of non-abelian simple groups are known,
see \cite[P. 302]{Gorenstein}.

\begin{lemma}\label{Double-Cov}{\rm(\cite[Proposition 2.6]{DFZ17})}
For $n\ge 7$, the covering group
$\ZZ_2.\A_n$ has no subgroup isomorphic to $\ZZ_2\times\A_{n-1}$.

\end{lemma}

A typical induction method for studying transitive graphs
is taking normal quotient graphs.
Suppose $\Ga$ is a $X$-vertex-transitive graph,
where $X\leq \Aut\Ga$ has an intransitive normal subgroup $N$.
Denote by $V\Ga_N:=\{\a^N\mid \a\in V\Ga\}$
the set of $N$-orbits in $V\Ga$. The {\it normal quotient graph}
$\Ga_N$ of $\Ga$ induced by $N$ is defined  with
vertex set $V\Ga_N$ and two vertices $B,C\in V\Ga_N$ are adjacent if and
only if some vertex in $B$ is adjacent in $\Ga$ to some vertex in
$C$. If $\Ga$ and $\Ga_N$ have the same valency,
then $\Ga$ is called a {\it normal cover}
(or {\it regular cover}) of $\Ga_N$.

The following proposition is a special case of \cite[Lemma 2.5]{Li}
which slightly improves a remarkable result of Praeger \cite[Theorem 4.1]{Praeger92}.

\begin{proposition}\label{praeger} Let $\Ga$ be a connected $X$-arc-transitive graph
of prime valence, with $X\le\Aut\Ga$, and let $N\lhd X$ have at least three orbits on
$V\Ga$. Then the following statements hold.
\begin{itemize}
\item[(1)] $N$ is semi-regular on $V\Ga$, $X/N\le\Aut\Ga_N$,
$\Ga_N$ is a connected $X/N$-arc-transitive graph, and $\Ga$ is a normal cover
of $\Ga_N$;
\item[(2)] $X_{\a}\cong(X/N)_{v}$ for $\a\in V\Ga$ and $v\in V\Ga_N$.
\end{itemize}
\end{proposition}

Let $G$ be a group, $g\in G$
and $H$ a subgroup of $G$.
Define the {\it coset graph} with vertex set $[G:H]$ (the set of cosets of $H$ in $G$),
and $Hx$ is adjacent to $Hy$ with $x,y\in G$ if and only if
$yx^{-1}\in HgH$. This coset graph is denoted by $\cos(G,H,g)$.
The following proposition follows from Sabidussi \cite{Sab64}.

\begin{proposition}\label{sab64}
Let $\Ga$ be a connected $G$-arc-transitive graph
of valency $k$, and $\a$ a vertex of $\Ga$.
Then $\Ga\cong\cos(G,G_{\a},g)$ with $g$ satisfying the following condition$:$

\vskip0.1in
{\noindent\bf Condition:} $g$ is a $2$-element of $G$,
$g^2\in G_{\a}$, $\l G_{\a},g\r=G$
and $k=[G_{\a}:G_{\a}\cap G_{\a}^g]$.

\vskip0.1in Conversely, if the element $g$ satisfies the above condition,
then $\cos(G,G_{\a},g)$ is a connected $G$-arc-transitive graph
of valency $k$.

\end{proposition}

Following the term in \cite{DFZ17},
the element $g$ satisfying the above condition
is called a {\it feasible element} to $G$ and $G_{\a}$.
The Magma \cite{Magma}
provides a powerful tool to compute out all the feasible elements $g$
giving rise to non-isomorphic graphs
(and so all the non-isomorphic graphs $\Ga$)
while the order of $G$ is not too large.
By direct computation, we have the following.

\begin{lemma}\label{Exm-1}
There is no connected $G$-arc-transitive $7$-valent graph
with $(G,G_{\a})$ lying in the following table.

\[\begin{array}{l|l|l|l|l} \hline
G &   \ZZ_2\times\A_7 & \ZZ_3\times\A_7 & \ZZ_4\times\A_7,~~ \ZZ_2^2\times\A_7
& \ZZ_9\times\A_7,~~ \ZZ_3^2\times\A_7   \\ \hline
G_{\a} & \F_{42}\times\ZZ_3 &\F_{42}\times\ZZ_2 & \F_{42}\times\ZZ_2,~~ \F_{42}\times\ZZ_3
& \F_{42}\times\ZZ_3,~~ \F_{42}\times\ZZ_6 \\ \hline
\end{array}\]

\end{lemma}

\section{A few technical lemmas}

In this section, we prove a few lemmas for proving Theorem~\ref{Thm-1}.
The first is a simple observation.

\begin{lemma}\label{simgp-mid-504}
Let $T$ be a non-abelian simple group with order dividing $504$.
Then $T\cong\PSL(2,7)$ and $\PSL(2,8)$,
with order $168$ and $504$ respectively.
\end{lemma}

Actually, since $504=2^3\cdot 3^2\cdot 7$,
$T$ is a $\{2,3,7\}$-group.
By \cite[Theorem III(2)]{Huppert-L},
the only $\{2,3,7\}$-nonabelian simple groups are
$\PSL(2,7)$, $\PSL(2,8)$ and $\PSU(3,3)$.
The lemma follows since $|\PSU(3,3)|=2^5\cdot 3^3\cdot 7$
does not divide $504$.

\vskip0.1in
The next depends on the classification of primitive
permutation groups with degree less than 1000,
obtained by Dixon and Mortimer.

\begin{lemma}\label{pri-less-252}
Let $S$ be a non-abelian simple group
with a non-abelian simple subgroup $T$ such that $|S:T|\div 252$
and $7\div |S|$.
Then either

\begin{itemize}
\item[(1)] $(S,T)=(\A_{n},\A_{n-1})$ where $n \geq 7$ divides $252$; or
\item[(2)] the couple $(S,T)$ is listed in the following table.
\end{itemize}

\begin{table}[ht]
\[\begin{array}{llll} \hline
S & \A_7 & \PSU(3,3) & \PSU(4,3) \\ \hline
T & \A_5 & \PSL(2,7) & \PSU(4,2) \\  \hline
\end{array}\]
\caption{`Sporadic' simple groups $S$ and $T$ with $|S:T|\div 252$}
\end{table}
\end{lemma}

\proof Let $H$ be a maximal subgroup of $S$ containing $T$,
and let $\Ome=[S:H]$.
Then $|\Ome|\div 252$
and $S$ acts faithfully and primitively on $\Ome$ by coset action,
that is,
$S$ can be viewed as a primitive permutation group on $\Ome$
with degree dividing $252$.
Inspecting the classification of primitive permutation groups with degree less than 1000
(see \cite[Appendix B]{DM}),
we conclude that either $(S,H)=(\A_{n},\A_{n-1})$ where $n \ge 7$ divides $252$,
or $(S,H)$ lies in the following table.

\begin{table}[ht]
\[\begin{array}{l|lllll} \hline
S & \A_7 & \PSL(3,4) & \A_8 & \PSp(6,2) &\A_9  \\
H & \S_5 & [2^4]:\A_5 & \S_6 & \PSU(4,2):2 & \S_7  \\
|S:H| & 21 & 21 & 28 & 28& 36 \\   \hline
S & \PSU(3,3) & \PSp(6,2) &   \PSL(6,2) & \PSp(6,2)& \A_9 \\
H & \PSL(2,7) & \S_8 & [2^5]:\PSL(5,2) & [2^5]:\S_6 & (\A_6 \times 3):2 \\
|S:H|  &  36  & 36 & 63 & 63 & 84\\  \hline

S & \A_9 & \PSp(6,2) & \PSU(4,3)& \\
H & (\A_5 \times \A_4):2 & (\A_5 \times \A_5):4& \PSU(4,2) &\\
|S:H| & 126 & 126& 126 & \\ \hline
\end{array}\]
\caption{Primitive permutation groups $S$ with degree dividing $252$ and insoluble stabilizer $H$}
\end{table}

Assume first $(S,H)=(\A_{n},\A_{n-1})$ with $n\geq 7$ dividing $252$.
If $H=T$, part (1) of Lemma~\ref{pri-less-252} holds.
Suppose $T$ is a proper subgroup of $H$.
Since $\A_{n-1}$ has no subgroup with index less that $n-1$,
$|H:T|\ge n-1$.
It follows $n(n-1)\le |S:H|\cdot |H:T|=|S:T|\le 252$,
implying $n\le 16$ because $17\cdot 16>252$.
Then since $n\ge 7$ and $n\div 252$, we further obtain $n=7,9,12$ or $14$.
Now, since $T<H\cong\A_{n-1}$ and $|\A_n:T|\div 252$,
by Atlas \cite{Atlas}, it is routing to check out
the only possibility is $(S,H,T)=(\A_7,\A_6,\A_5)$.

Assume now $(S,H)$ lies in Table 2.
Noting that $|S:H|$ is listed in Table 2
and $|H:T|$ divides ${252\over |S:H|}$.
If $(S,H)=(\A_7,\S_5)$ or $(\PSU(3,3),\PSL(2,7))$,
one easily has $T=\A_5$ or $\PSL(2,7)$  respectively;
and if $(S,H)=(\PSU(4,3),\PSU(4,2))$,
then $|H:T|\div 2$ because $|S:H|=126$,
implying $T=H=\PSU(4,2)$,
which give rise to examples in Table 1.
Suppose $(S,H)=(\PSp(6,2),\PSU(4,2))$.
Then $|S:H|=28$,
and so $|H:T|\div 9$,
by \cite{Atlas}, no such a group $T$ exists.
For the remaining cases in Table 2,
one may similarly check out no group $T$ exists.\qed

As usual, for a group $R$,
we denote by $R_p$ and $R_{\pi}$ a Sylow $p$-subgroup
and a Hall $\pi$-subgroup of $R$,
where $p$ is a prime dividing $|R|$ and $\pi$ is a set of some primes dividing $|R|$.

\begin{lemma}\label{252}
Let $R$ be a soluble $\{2,3,7\}$-group
such that $|R|\div 2^2\cdot 3^2\cdot7$. Then
either $R_{2,7}\lhd R$ or $R_{3,7}\lhd R$.
\end{lemma}

\proof
Since $R$ is soluble,
$R$ has a Sylow system $R_2,R_3,R_7$,
that is, $R_iR_j$ is a subgroup of $R$ for distinct $i,j\in\{2,3,7\}$.
Since $|R_2|\div 4$ and $|R_3|\div 9$, by Sylow theorem,
both $R_2R_7$ and $R_3R_7$ has a unique Sylow $7$-subgroup $R_7$,
that is, both $R_2$ and $R_3$
normalize $R_7$.
It follows $R_7\lhd R$.

Set $H=R_2R_3$.
Assume $|H|=36$. Let $F$ be the {\it Fitting subgroup} of $H$,
that is, the largest nilpotent normal subgroup of $H$.
Clearly, $F$ is abelian.
Then since $H$ is soluble, by \cite[P. 30, Corollary]{Suzuki},
$F\ne 1$ and $C_R(F)=F$.
If $|F|\div 6$, then $F\le\ZZ_6$
and Lemma~\ref{N/C} implies $H/F\lesssim \Aut(F)\le\ZZ_3$,
which is a contradiction because $|H/F|\ge 6$.
Thus, $|F|$ is divisible by $4$ or $9$.
For the former case, $R_2=F_2\lhd H$,
and so $R_{2,7}\lhd R$.
For the latter case, $R_3=F_3\lhd H$
and $R_{3,7}\lhd R$.
If $|H|\div 12$ or $18$, one easily shows
that one of $R_2$ and $R_3$ is normal in $H$,
the lemma follows.\qed

The {\it soluble radical} of a group is its largest soluble normal subgroup.

\begin{lemma}\label{small-T}
Let $\Ga$ be a connected arc-transitive $7$-valent Cayley graph on $T$,
where $T\le\A_7$ is a non-abelian  simple group,
and let $R$ be the soluble radical of $\Aut\Ga$.
If $R\ne 1$ is semi-regular on $V\Ga$ and $(\Aut\Ga)_{\a}$
with $\a\in V\Ga$ is soluble, then $RT=R\times T$.
\end{lemma}

\proof By assumption, $T=\A_5,\A_6,\A_7$ or $\PSL(2,7)$.
Let $Y=RT$. Since $R\cap T\lhd T$ is soluble,
$R\cap T=1$ and $Y=R:T$.
Since $|V\Ga||Y_{\a}|=|Y|=|R||T|=|R||V\Ga|$,
by Lemma~\ref{Stabilizer}, one has $|R|=|Y_{\a}|$ divides $252=2^2\cdot3^2\cdot7$.
If $|R|$ is square-free or $|R|$ has at most two prime divisors,
one easily shows that $\Aut(R)$ is soluble.
Then by Lemma~\ref{N/C}, $Y/C_Y(R)\lesssim\Aut(R)$,
implying $T\le C_Y(R)$,
hence $Y=R\times T$, as required.

The remaining cases are $|R|=2^2\cdot 3\cdot 7=84$, $2\cdot 3^2\cdot 7=126$
and $2^2\cdot 3^2\cdot 7=252$.
In these cases, $|Y_{\a}|=|R|$ is divisible by $7$,
so $\Ga$ is $Y$-arc-transitive,
and since $R$ is semi-regular on $V\Ga$, $|R|$ divides $|V\Ga|=|T|$,
in particular $T\ne\A_5$ and $\A_6$.
Since a normal Hall subgroup is a characteristic subgroup,
by Lemma~\ref{252}, one of $R_{2,7}$ and $R_{3,7}$ is normal in $Y$.
Clearly, both $R_{2,7}$ and $R_{3,7}$ have at least three orbits on $V\Ga$.

Suppose $T=\PSL(2,7)$. Then $|V\Ga|=|\PSL(2,7)|=168$ and $|R|=84$.
If $R_{2,7}\lhd Y$, by Proposition~\ref{praeger}, the normal quotient graph
$\Ga_{R_{2,7}}$ is an arc-transitive 7-valent graph of order ${|V\Ga|\over |R_{2,7}|}=6$, a contradiction.
If $R_{3,7}\lhd Y$, again by Proposition~\ref{praeger}, $Y/R_{3,7}\le\Aut(\Ga_{R_{3,7}})$,
and $\Ga_{R_{3,7}}$ is an $Y/R_{3,7}$-arc-transitive 7-valent graph
of order ${|V\Ga|\over |R_{3,7}|}=8$,
hence $\Ga_{R_{3,7}}\cong\K_8$ is a complete graph;
however, as $|R_2|=4$, one has $Y/R_{3,7}\cong R_{2}:\PSL(2,7)\cong R_2\times\PSL(2,7)$,
and by Atlas \cite{Atlas}, $\Aut(\K_8)\cong\S_8$ has no such a subgroup,
also a contradiction.

Suppose now $T=\A_7$.
Then $|V\Ga|=|\A_7|=2520$, and $|R|=84,126$ or $252$.

Assume $|R|=84$. Then $|R_{2,7}|=28$ and $|R_{3,7}|=21$,
and as $|Y_{\a}|=|R|$, we have $Y_{\a}\cong\F_{42}\times\ZZ_2$ by Lemma~\ref{Stabilizer}.
If $R_{2,7}\lhd Y$,
Proposition~\ref{praeger} implies that
$\Ga_{R_{2,7}}$ is a connected $Y/R_{2,7}$-arc-transitive 7-valent graph,
and $(Y/R_{2,7})_v\cong Y_{\a}\cong \F_{42}\times\ZZ_2$ for $v\in V(\Ga_{R_{2,7}})$;
also, since $Y=R:T$ and $R/R_{2,7}\cong\ZZ_3$,
we have $Y/R_{2,7}\cong\ZZ_3:\A_7\cong\ZZ_3\times\A_7$.
By Lemma~\ref{Exm-1}, no such graph $\Ga_{R_{2,7}}$ exists,
a contradiction.
Similarly, if $R_{3,7}\lhd Y$,
then $\Ga_{R_{3,7}}$ is a connected $Y/R_{3,7}$-arc-transitive 7-valent graph,
$(Y/R_{3,7})_v\cong \F_{42}\times\ZZ_2$
and $Y/R_{3,7}\cong R_2:\A_7
\cong\ZZ_4\times\A_7$ or $\ZZ_2^2\times\A_7$,
by Lemma~\ref{Exm-1}, no such graph $\Ga_{R_{3,7}}$ exists,
also a contradiction.

For the cases where $|R|=126$ and $252$,
by Lemma~\ref{Stabilizer}, one has $Y_{\a}\cong\F_{42}\times\ZZ_3$
and $\F_{42}\times\ZZ_6$ respectively.
Then, with similar discussion as above,
we have that $\Ga_{R_{i,7}}$ is a
connected $Y/R_{i,7}$-arc-transitive 7-valent graph for $i=2$ or $3$,
with the tuple $(Y/R_{i,7},(Y/R_{i,7})_v)$ lying in the table of Lemma~\ref{Exm-1},
it is a contradiction.\qed

We next prove a more general assertion.

\begin{lemma}\label{R-not-1}
Let $\Ga$ be a connected $7$-valent arc-transitive Cayley graph on
a non-abelian simple group $T$,
and let $R$ be the soluble radical of $\Aut\Ga$.
If $R\ne 1$ is semi-regular on $V\Ga$ and $(\Aut\Ga)_{\a}$
is soluble for $\a\in V\Ga$, then $RT=R\times T$.
\end{lemma}

\proof Set $Y=RT$. Then $Y=R:T\le\Aut\Ga$
and $|R|=|Y_{\a}|$ divides $2^2\cdot 3^2\cdot 7$.
If $|R|$ has at most two prime divisors,
by the arguments in the proof of Lemma~\ref{small-T}, one has
$Y=R\times T$.
If $T\le\A_7$, by Lemma~\ref{small-T}, $Y=R\times T$.

Thus assume that in the following, $|R|\ne 7$, $7\div |R|$ and $T{\not\le}\A_7$.
Then $|Y_{\a}|=|R|$ is divisible by $7$
and $\Ga$ is $Y$-arc-transitive.
Since $R$ is soluble and $|R|\div 2^2\cdot 3^2\cdot 7$,
$R$ has a nontrivial Hall $\{2,3\}$-group, say $H$.
Let $\Ome=\{ H^r\div r\in R\}$, the set of conjugate subgroups of $H$ under $R$.
Then $|\Ome|=|R:N_R(H)|=1$ or $7$,
and the conjugate action of $Y$ on $\Ome$ is transitive.
If $T$ acts faithfully on $\Ome$, then $|\Ome|=7$
and $T\le \A_7$, which is not the case.

Assume $T$ acts unfaithfully on $\Ome$.
Since $T$ is non-abelian simple, the conjugate action of $T$ on $\Ome$ is trivial, that is,
$T$ normalizes $H$, hence $TH\le Y$.
Set $Z=TH$ and let $\Del=[Y:Z]$.
Since $H\cap T\le R\cap T=1$, we have $|Z|=|H||T|$
and $|\Del|=|Y:Z|=7$.
Let $Y_{(\Del)}$ be the kernel of $Y$ acting on $\Del$ by coset action.
Then $Y_{(\Del)}=\cap_{y\in Y}Z^y$ is the largest normal subgroup of $Y$ contained in $Z$,
and $Y/Y_{(\Del)}\le\S_7$.
It follows $T\le Y_{(\Del)}$ because $T{\not\le}\A_7$.
Since $|H||T|=|Z|=|T||Z_{\a}|$,
$Z_{\a}$ is a $\{2,3\}$-group as $H$ is.
If $(Y_{(\Del)})_{\a}\ne 1$,
since $\Ga$ is a connected $Y$-arc-transitive graph of valency $7$
and $1\ne (Y_{(\Del)})_{\a}\lhd Y_{\a}$, we have
$1\ne (Y_{(\Del)})_{\a}^{\Ga(\a)}\lhd Y_{\a}^{\Ga(\a)}$,
implying $7\div |(Y_{(\Del)})_{\a}|$.
Hence $7\div |Z_{\a}|$ as $Y_{(\Del)}\subseteq Z$, a contradiction.
Thus $(Y_{(\Del)})_{\a}=1$. Now,
since $T\le Y_{(\Del)}$, by Lemma~\ref{Frattini}, $T=T(Y_{(\Del)})_{\a}=Y_{(\Del)}$
is normal in $Y$, hence $Y=R\times T$.\qed

\section{Proof of Theorem~\ref{Thm-1}}

In this section, we will complete the proof of Theorem~\ref{Thm-1}.
Recall the {\it socle} of a group is the product of its all minimal normal subgroups.

\begin{lemma}\label{R=1}
Let $\Ga$ be a connected $X$-arc-transitive $T$-vertex-transitive $7$-valent graph,
where $T\le X\le\Aut\Ga$ is a non-abelian simple group.
Suppose $X_{\a}$ is soluble for $\a\in V\Ga$ and the soluble radical of $X$ is trivial.
Then either $T\lhd X$ or $X$ is almost simple with socle $S$ an overgroup of $T$
such that $(S,T)=(\A_n,\A_{n-1})$, where $n\ge 7$ is a divisor of $252$.
\end{lemma}

\proof Suppose $T$ is not normal in $X$.
Let $N$ be a minimal normal subgroup of $X$.
Since the trivial soluble radical of $X$ is trivial,
$N=S^d$ with $S$ a non-abelian simple group and $d\ge 1$.

Let $Y=NT$. Then $Y\le X$, and  $|Y|={|N||T|\over |N\cap T|}$
divides $|X|=|TX_{\a}|={|T||X_{\a}|\over|T\cap X_{\a}|}=|T||X_{\a}:T_{\a}|$.
It follows ${|N|\over |N\cap T|}$ divides
$|X_{\a}|$, and hence divides $252$ by Lemma~\ref{Stabilizer}.
Since $N\cap T\lhd T$,
$N\cap T=1$ or $T$.
If $N\cap T=1$, then $|N|=|S|^d$ divides $252$,
which is a contradiction by Lemma~\ref{simgp-mid-504}.
Hence $T\le N$, $Y=N$ and $|N:T|\div 252$.
Write $N=S_1\times S_2\times\cdots\times S_d$ with each $S_i\cong S$.
Assume $d\ge 2$. If $S_1\cap T=1$,
then $S_1=S_1/(S_1\cap T)\cong S_1T/T\le N/T$,
one has $|S_1|\div 252$, by
Lemma~\ref{simgp-mid-504}, it is a contradiction.
If $S_1\cap T\ne 1$,
since $S_1\cap T\lhd T$, we have $T\le S_1$
and $|S_2|\div |N:S_1|\div |N:T|$.
Hence $|S_2|\div 252$, also a contradiction by Lemma~\ref{simgp-mid-504}.
Thus $d=1$ and $S=N\ge T$.
Further, if $X$ has another minimal normal subgroup $M$,
then the above discussion may imply that $M\ge T$ is simple.
It follows $MN=M\times N\le X$
and $|M|\div 252$, by Lemma~\ref{simgp-mid-504},
it is a contradiction.

Therefore, $X$ is almost simple with socle $S\ge T$.
Since $T$ is not normal in $X$, $S>T$ and $S_{\a}\ne 1$.
Since $\Ga$ is connected and $1\ne S_{\a}\lhd X_{\a}$,
one has $1\ne S_{\a}^{\Ga(\a)}\lhd X_{\a}^{\Ga(\a)}$.
It follows $7\div |S_{\a}|$ and $\Ga$ is $S$-arc-transitive.
Since $S=TS_{\a}$, we have $|S|=|TS_{\a}|={|T||S_{\a}|\over|T_{\a}|}$,
so $|S_\a:T_\a|=|S:T|$. In particular, $|S:T|\div |S_{\a}|\div 252$.
Hence $(S,T)$ satisfies Lemma \ref{pri-less-252}.

If $(S,T)=(\A_7,\A_5)$, then $|S:T|=42$ divides $|S_\a|$.
By Lemma~\ref{Stabilizer}, $S_\a\cong\F_{42}\times \ZZ_m$ with $m=1,2,3$ or $6$.
A direct computation by Magma \cite{Magma} shows
that no feasible element exists for $S$ and $S_{\a}$,
that is, no graph $\Ga$ exists, a contradiction.
For the cases where $(S,T)=(\PSU(3,3),\PSL(2,7))$ and $(\PSU(4,3),\PSU(4,2))$,
one may similarly draw a contradiction.

Hence $(S,T)=(\A_n,\A_{n-1})$
with $n\ge 7$ a divisor of $252$,
the lemma follows.\qed

\begin{lemma}\label{R-not-1-case}
Let $\Ga$ be a connected arc-transitive $7$-valent Cayley graph on $T$,
and let $R$ be the soluble radical of $\Aut\Ga$, where
$T$ is a non-abelian simple group.
If $R\ne 1$ and $\A_{\a}$ is soluble for $\a\in V\Ga$,
then either $T\lhd \Aut\Ga$ or
$\Aut\Ga$ has a non-abelian normal simple subgroup $S>T$
such that $(S,T)=(\A_n,\A_{n-1})$, with $n\ge 7$ a divisor of $252$.
\end{lemma}

\proof Set $\A=\Aut\Ga$ and suppose $T$ is not normal in $\A$.
Let $Y=RT$. Then $Y=R:T$, and by Lemma~\ref{Stabilizer}, $|R|=|Y_{\a}|$ divides $252$.

Assume first $R$ is transitive on $V\Ga$.
Then $|V\Ga||Y_{\a}|=|Y|=|R||T|=|V\Ga||R_{\a}||T|$,
we obtain $|T|\div 252$, by Lemma~\ref{simgp-mid-504},
it is a contradiction.

Assume now  $R$ has exactly two orbits on $V\Ga$.
Then $|V\Ga||Y_{\a}|={1\over 2}|V\Ga||R_{\a}||T|$,
so $|Y_{\a}|={1\over 2}|R_{\a}||T|$
and $|T|\div 504$,
by Lemma~\ref{simgp-mid-504}, $T=\PSL(2,7)$ or $\PSL(2,8)$.
In particular, $7\div |T|$.
It follows
$7\div |Y_{\a}|$ and $\Ga$ is $Y$-arc-transitive.
If $R_{\a}\ne 1$, since $R\lhd Y$,
we have $1\ne R_{\a}^{\Ga(\a)}\lhd Y_{\a}^{\Ga(\a)}$
and hence $7\div |R_{\a}|$, thus $7^2\div |Y_{\a}|$,
it is a contradiction by Lemma~\ref{Stabilizer}.

Consider the case where $R_{\a}=1$. If $T=\PSL(2,7)$, then $|R|={1\over 2}|T|=84$,
by the arguments in the proof of Lemma~\ref{small-T},
no graph $\Ga$ exists in the case.
Suppose $(|R|,T)=(252,\PSL(2,8))$.
By Lemma~\ref{252}, one of $R_{2,7}$ and $R_{3,7}$ is normal in $Y$.
If $R_{2,7}\lhd Y$, then the normal quotient graph
$\Ga_{R_{2,7}}$ is a connected arc-transitive 7-valent graph of order $18$,
by \cite{MR90},
no such a graph exists, a contradiction;
if $R_{3,7}\lhd Y$, then $\Ga_{R_{3,7}}$ is an $Y/R_{3,7}$-arc-transitive 7-valent graph of order $8$,
hence $\Ga_{R_{3,7}}\cong\K_8$ is a complete graph,
however, as $|R_2|=4$, by Atlas \cite{Atlas}, $\Aut(\K_8)\cong\S_8$ has no subgroup
isomorphic to $Y/R_{3,7}\cong R_{2}:\PSL(2,8)$, also a contradiction.

Thus, assume in the following, $R$ has at least three orbits on $V\Ga$.
By Proposition~\ref{praeger}, $R$ is semi-regular on $V\Ga$
and so $Y=R\times T$ by Lemma~\ref{R-not-1},
and $\Ga_R$ is $\A/R$-arc-transitive and $TR/R$-vertex-transitive $7$-valent graph of order ${|V\Ga|\over |R|}$.
Since $Y=R\times T$, $T$ is a characteristic subgroup of $Y$,
then as $T$ is not normal in $\A$, we conclude $Y$ is not normal in $\A$,
and so $Y/R$ is not normal in $\A/R$.
Noting that the soluble radical of $\A/R$ is trivial, and $Y/R\cong T$,
by Lemma~\ref{R=1}, we have
$\A/R$ is almost simple (say with socle $X/R$)
such that $S\cong X/R> Y/R\cong T$ and $(X/R,Y/R)=(\A_{n},\A_{n-1})$
with $n\ge 7$ a divisor of $252$.

Let $C=C_X(R)$, the centralizer of $R$ in $X$. Then $C\lhd X$.
Since $Y=R\times T<X$, $C\ge T$ is insoluble and $C\cap R=Z(R)\le Z(C)$.
Since $C/(C\cap R) \cong CR/R \lhd X/R \cong S$, we obtain $C/(C\cap R) \cong S$.
Thus $C=(C\cap R).S$ is a central extension.
If $C\cap R <Z(C)$, then $1 \ne Z(C)/(C\cap R)\lhd C/(C\cap R) \cong S$, implying $Z(C)=C$, a contradiction.
Hence $C\cap R=Z(C)$ and $X=RC$.

Now, since $C'\cap Z(C)\le Z(C')$, we have $Z(C')/(C'\cap Z(C))\lhd C'/(C'\cap Z(C))\cong C'Z(C)/Z(C)=(C/Z(C))'\cong S'=S$.
It follows $C'\cap Z(C)=Z(C')$, $C=C'Z(C)$ and $C'=Z(C').S$.
Hence $C'=(C'Z(C))'=C''$,
that is, $C'$ is a covering group of $S$.
Hence $C'\cap Z(C)=Z(C')\le \Mult(S)$.
Since $T<C$ and $C/C'$ is abelian,
$T \leq C'$.

Assume first $S=\A_{n}$ with $n>7$.
By \cite[Theorem 5.14]{KL90}, $\Mult(S)\cong\ZZ_2$,
thus $Z(C')=1$ or $\ZZ_2$.
If $Z(C')\cong\ZZ_2$,  then $C'=2.S$. Since $T\cap Z(C')=1$,
$T Z(C')\cong \A_{n-1}\times\ZZ_2$ is a subgroup of $C'$;
however, by Lemma~\ref{Double-Cov}, $C'\cong 2.\A_{n} $ has no subgroup isomorphic to $\A_{n-1} \times \ZZ_2$,
it is a contradiction.
Thus $Z(C')=1$. Now, $C'\cong S$ and $C=Z(C)C'=Z(C)\times C'$.
Recall that $X/R\cong S$ and $X=RC$ has a subgroup $RC'\cong R\times S$,
we conclude $X=RC'=R \times C'\lhd \A$.
It follows $C'\cong S$ is a characteristic subgroup of $X$. Hence $C'\lhd A$,
with $(C',T)=(S,T)=(\A_n,\A_{n-1})$, the lemma is true.

Assume now $S=\A_7$. Then $Z(C')\le\Mult(\A_7)\cong\ZZ_6$
by \cite[Theorem 5.14]{KL90}.
With above discussions, one may exclude the case where $Z(C')\cong\ZZ_2$,
and if $Z(C')=1$, the lemma is true with $(S,T)=(\A_7,\A_6)$.
If $Z(C')\cong\ZZ_6$, then the covering group $C'\cong\ZZ_6.\A_7$ has a subgroup
$T Z(C')\cong \A_6\times\ZZ_6$.
It follows $C'$ has a quotient group which is isomorphic to a covering group
$\ZZ_2.\A_7$ and contains a subgroup $\A_6\times\ZZ_2$, a contradiction by Lemma~\ref{Double-Cov}.
If $Z(C')\cong\ZZ_3$, then $C'\cong\ZZ_3.\A_7$ has a subgroup
$G Z(C')\cong \A_{6}\times\ZZ_3$,
but a direct computation by Magma \cite{Magma} shows that the covering roup
$\ZZ_3.\A_7$ has no subgroup isomorphic to $\A_{6}\times\ZZ_3$,
it is a contradiction.\qed

Now, we are ready to prove Theorem~\ref{Thm-1}.

\vskip0.1in
{\noindent\bf Proof of Theorem~\ref{Thm-1}.}
Set $\A=\Aut\Ga$.
Then $\Ga$ is $\A$-arc-transitive.
By Lemmas~\ref{R=1} and ~\ref{R-not-1-case},
either $T\lhd\A$ or $\A$ has a normal subgroup $S>T$
such that $(S,T)=(\A_n,\A_{n-1})$,
where $n\ge 7$ is a divisor of $252$.
For the former case,
part (1) of Theorem~\ref{Thm-1} holds.

Consider the latter case.
Since $T$ is regular on $V\Ga$, $|S_{\a}|=n$ for $\a\in V\Ga$.
Then since $S\lhd \A$,
we have $1\ne S_{\a}\lhd\A_{\a}$,
implying $1\ne S_{\a}^{\Ga(\a)}\lhd\A_{\a}^{\Ga(\a)}$.
It follows that $S_{\a}^{\Ga(\a)}$ is transitive,
and $\Ga$ is $S$-arc-transitive.
In particular, $7\div n$. Since $n\div 252$,
$n\in\{7,14,21,28,42,63,84,126,252\}$.

For the cases $n=14,42$ and $126$, let $\Ome=[S:T]$,
and consider the coset action of $S$ on $\Ome$.
Since $|\Ome|=n$ and $S=\A_n$ has a unique conjugate class
of subgroups isomorphic to $\A_{n-1}$,
this action is permutation equivalent to the natural action of $\A_n$ on $n$ points.
Since $S=TS_{\a}$ and $|S_{\a}|=n$,
$S_{\a}$ acts regularly on $\Ome$.
However, since $S_{\a}\le\A_n$ is of order twice an odd integer,
any involutions in $S_{\a}$ will fix some point in $\Ome$,
which is a contradiction.
If $n=28$, then $S=\A_{28}$, and by Lemma~\ref{Stabilizer}, $S_{\a}\cong\D_{14}\times\ZZ_2$.
A direct computation by Magma \cite{Magma} shows
there is no feasible element to $\A_{28}$ and $\D_{14}\times\ZZ_2$,
thus no graph $\Ga$ exists in the case.
Similarly, if $n=252$, then $S=\A_{252}$ and $S_{\a}\cong\F_{42}\times\ZZ_6$,
by Magma \cite{Magma}, there exists no feasible element to $\A_{252}$ and $\F_{42}\times\ZZ_6$,
hence no graph $\Ga$ exists in the case.

Therefore, $n=7,21,63$ or $84$. By \cite[Theorem 1.3]{FMW11},
for each prime $p>5$, there is a connected $p$-valent non-normal
$\A_p$-arc-transitive Cayley graph on $\A_{p-1}$,
so $\Ga$ exists for the case $n=7$;
and if $n=21,63$ and $84$, by Examples~\ref{21}-\ref{84} below,
there exist connected
$\A_n$-arc-transitive $7$-valent non-normal
Cayley graphs on $\A_{n-1}$,
the last statement of Theorem~\ref{Thm-1} is true.
Finally, if $n=7,21$ or $63$, by Lemma~\ref{Stabilizer},
$S_{\a}\cong\ZZ_7,\F_{21}$
or $\F_{21}\times\ZZ_3$ respectively,
and $\Ga$ is $(S,1)$-transitive, as in part (2)(i) of Theorem~\ref{Thm-1}.
If $n=84$,  then Lemma~\ref{Stabilizer} implies $S_{\a}=\F_{42}\times\ZZ_2$,
and $\Ga$ is $(\A_{84},2)$-transitive,
part (2)(ii) of Theorem~\ref{Thm-1} holds.\qed

\section{Examples and the full automorphism groups}

In this section, we give some specific examples of graphs
satisfying parts (2)(i) and (2)(ii) of Theorem~\ref{Thm-1}.

\begin{example}\rm \label{21}
Let $S=\A_{\lbrace1,2,...,21\rbrace} \cong \A_{21}$, $T=\A_{\lbrace2,3,...,21\rbrace} \cong \A_{20}$ and
\begin{itemize}
\item[$x$]=(1, 2, 4)(3, 6, 10)(5, 9, 14)(7, 12, 8)(11, 17, 20)(13, 19, 16)(15, 21, 18);
\item[$y$]=(1, 3, 7, 13, 14, 20, 21)(2, 5, 6, 11, 12, 18, 19)(4, 8, 9, 15, 10, 16, 17);
\item[$g$]=(7, 11)(8, 20)(12, 17)(13, 15)(16, 18)(19, 21).
\end{itemize}

Let $H=\l x,y\r$ and $\Ga=\cos(S,H,g)$.

\vskip0.1in By Magma \cite{Magma}, $H=\l y\r: \l x\r \cong\F_{21}$, $\l H,g\r=S$
and $|H:H \cap H^g|=7$.
By Proposition~\ref{sab64}, $\Ga$ is a connected $\A_{21}$-arc-transitive $7$-valent graph.
Noting that $\l y\r$ has $3$ orbits on $\{1,2,\dots,21\}$,
$\l x\r$ permutates cyclically these three orbits and $|H|=21$,
we conclude $H$ is regular on $\lbrace1,2,...,21\rbrace$.
Since the vertex stabilizer $S_1=T$, by Lemma~\ref{Frattini},
$S=HT$. It follows that $T$ is regular on $V\Ga=[S:H]$,
that is, $\Ga$ is a Cayley graph on $\A_{20}$.
Finally, since $\A_{21}\cong S\le\Aut\Ga$, $\Ga$ is non-normal.
\end{example}

\begin{example}\rm \label{63}
Let $S=\A_{\lbrace1,2,...,63\rbrace} \cong \A_{63}$, $T=\A_{\lbrace2,3,...,63\rbrace} \cong \A_{62}$ and
\begin{itemize}
\item[$x$]=(1, 2, 5)(3, 6, 12)(4, 7, 13)(8, 14, 24)(9, 15, 25)(10, 16, 26)(11, 17, 27)(18, 28, 40)(19, 29, 41)(20, 30, 42)(21,31, 43)(22, 32, 44)(23, 33, 45)(34, 46, 54)(35, 47, 55)(36, 48, 56)(37, 49, 57)(38, 50, 58)(39, 51, 59)(52, 60, 62)(53, 61, 63);
\item[$y$]=(1, 3, 8)(2, 6, 14)(4, 10, 20)(5, 12, 24)(7, 16, 30)(9, 19, 34)(11, 22, 18)(13, 26, 42)(15, 29, 46)(17, 32, 28)
	(21,37, 52)(23, 39, 36)(25, 41, 54)(27, 44, 40)(31, 49, 60)(33, 51, 48)(35, 53, 38)(43, 57, 62)(45, 59, 56)(47, 61, 50)(55, 63, 58);
\item[$z$]=   (1, 4, 11, 23, 34, 52, 53)(2, 7, 17, 33, 46, 60, 61)(3, 9, 10, 21, 22, 38, 39)(5, 13, 27, 45, 54, 62, 63)(6, 15, 16,31, 32, 50, 51)
	(8, 18, 19, 35, 20, 36, 37)(12, 25, 26, 43, 44, 58, 59)(14, 28, 29, 47, 30, 48, 49)(24, 40, 41,55, 42, 56, 57);
\item[$g$]=   (2, 3)(5, 8)(7, 10)(12, 14)(13, 20)(15, 19)(17, 22)(18, 27)(23, 35)(25, 34)(26, 30)(28, 44)(31, 37)(33, 53)(36, 55)(38,
    45)(39, 47)(41, 46)(43, 52)(48, 63)(50, 59)(51, 61)(56, 58)(57, 60).
\end{itemize}

Let $H=\l x,y,z\r$ and let $\Ga=\cos(S,H,g)$.

\vskip0.1in By Magma \cite{Magma}, $H=\l x\r\times (\l z\r: \l y\r)\cong\ZZ_3\times\F_{21}$, $\l H,g\r=S$
and $|H:H \cap H^g|=7$.
By Proposition~\ref{sab64}, $\Ga$ is a connected $\A_{63}$-arc-transitive $7$-valent graph.
Also, it is easy to show that $H$ is regular on $\lbrace1,2,...,63\rbrace$.
Since the vertex stabilizer $S_1=T$, by Lemma~\ref{Frattini},
$S=HT$, hence $T$ is regular on $V\Ga=[S:H]$,
that is, $\Ga$ is a Cayley graph on $\A_{62}$.
Finally, since $T=\A_{62}$ is not normal in $S=\A_{63}$,
$\Ga$ is non-normal.
\end{example}

\begin{example}\rm \label{84}
Let $S=\A_{\lbrace1,2,...,84\rbrace} \cong \A_{84}$, $T=\A_{\lbrace2,3,...,84\rbrace} \cong \A_{83}$ and
\begin{itemize}
\item[$x$]=(1, 3)(2, 7)(4, 13)(5, 11)(6, 16)(8, 22)(9, 20)(10, 25)(12, 29)(14, 34)(15, 32)(17, 39)(18, 38)(19, 42)(21, 45)(23, 48)(24, 50)(26, 55)(27, 53)(28, 58)(30, 54)(31, 61)(33, 51)(35, 66)(36, 65)(37, 62)(40, 69)(41, 59)(43, 71)(44, 56)(46, 49)(47, 73)(52, 76)(57, 79)(60, 81)(63, 78)(64, 77)(67, 84)(68, 83)(70, 82)(72, 80)(74, 75);
\item[$y$]=(1, 63, 21, 5, 30, 47)(2, 74, 14, 9, 46, 36)(3, 77, 41, 11, 51, 70)(4, 76, 62, 15, 50, 83)(6, 67, 12, 18, 35, 31)(7, 81, 26, 20, 58, 57)(8, 82, 25, 23, 59, 53)(10, 61, 39, 27, 29, 69)(13, 43, 45, 32, 19, 73)(16, 80, 24, 38, 56, 52)(17, 79, 42, 40, 55, 71)(22, 68, 34, 48, 37, 65)(28, 72, 33, 60, 44, 64)(49, 84, 54, 75, 66, 78);
\item[$g$]=(1, 6)(3, 18)(4, 63)(5, 16)(7, 9)(8, 82)(10, 77)(11, 38)(12, 79)(13, 30)(14, 76)(15, 78)(17, 74)(19, 73)(21, 71)(22, 59)(23, 70)(24, 65)(25, 51)(26, 61)(27, 64)(28, 83)(29, 55)(31, 57)(32, 54)(33, 53)(34, 50)(35, 56)(36, 52)(37, 81)(39, 46)(40, 75)(41, 48)(42, 45)(43, 47)(44, 67)(49, 69)(58, 62)(60, 68)(66, 80).
\end{itemize}

Let $H=\l x,y\r$ and let $\Ga=\cos(S,H,g)$.

\vskip0.1in By Magma \cite{Magma}, $H\cong \ZZ_2 \times \F_{42}$, $\l H,g\r=S$
and $|H:H \cap H^g|=7$.
Hence Proposition~\ref{sab64} implies that $\Ga$ is a connected $\A_{84}$-arc-transitive $7$-valent graph.
Also, with similar discussion as above, one
has that $H$ is regular on $\lbrace1,2,...,84\rbrace$,
and $\Ga$ is a non-normal Cayley graph on $T=\A_{83}$.
\end{example}

At the end of this paper, we determine the full automorphism group
of the graph in Example~\ref{21}. One may similarly
approach the full automorphism groups
of the graphs in Examples~\ref{63} and \ref{84}.
Recall a transitive permutation group is called {\it quasiprimitive}
if each of its minimal normal subgroups is transitive.

\begin{lemma}\label{rem-d10-s10}
Let $\Ga=\Cos(S,H,g)$ be as in Example~$\ref{21}$. Then $\Aut\Ga\cong\A_{21}$
and $\Ga$ is $1$-transitive.
\end{lemma}

\proof Recall that $\A_{20}\cong T<S\cong\A_{21}$
and $\Ga$ is a connected $S$-arc-transitive $7$-valent Cayley graph on $T$.
Let $\A=\Aut\Ga$ and $\a\in V\Ga$.
By \cite[Theorem 1.1]{GLH16}, $|\A_{\a}|\div2^{24}\cdot3^4\cdot5^2\cdot7$.

Assume $\A$ is not quasiprimitive on $V\Ga$.
Then $\A$ has an intransitive minimal normal subgroup $N$.
Set $F=NS$.
Since $S$ is non-abelian simple and $N\cap S\lhd S$, $N\cap S=1$ or $S$.
If $N\cap S=S$, then $N$ is transitive on $V\Ga$, a contradiction.
Suppose $N\cap S=1$. Then $F=N:S$ and $|N|=|F:S|$ divides $|\A:S|$. Since
$|V\Ga|=|\A:\A_{\a}|=|S:S_{\a}|$, we have
$|\A:S|=|\A_\a:S_\a|$ divides $2^{24}\cdot3^3\cdot5^2$,
so is $|N|$.
Since $|V\Ga|=|T|=|\A_{20}|$, $N$ has at least three orbits on $V\Ga$.
By Proposition~\ref{praeger}, $N$ is semi-regular on $V\Ga$,
and so $|N|$ divides $|V\Ga|=|\A_{20}|=2^{17}\cdot3^8\cdot5^4\cdot7^2\cdot11\cdot13\cdot17\cdot19$.

If $N$ is insoluble, since $|N|\div 2^{24}\cdot3^3\cdot5^2$,
and $\A_5,\A_6$ and $\PSp(4,3)$ are the only $\{2,3,5\}$-simple groups (see \cite[TABLE 1]{Huppert-L}),
by checking the orders, we conclude $N\cong\A_5$, $\A_5^2$ or $\A_6$.
Then since $|N||\A_{21}|=|N||S|=|F|=|V\Ga||F_{\a}|=|\A_{20}||F_{\a}|$,
we have $|F_{\a}|=2^2\cdot3^2\cdot5\cdot7$, $2^4\cdot3^3\cdot5^2\cdot7$ or $2^3\cdot3^3\cdot5\cdot7$.
However, checking the orders of the stabilizers of connected 7-valent arc-transitive graphs
given in Lemma~\ref{Stabilizer},
it is a contradiction

Suppose $N$ is soluble. Noting that $|N|\div |\A_\a:S_{\a}|$,
$|\A_\a:S_{\a}|\div 2^{24}\cdot3^3\cdot5^2$
and $2^{18}$ does not divide $|N|$,
checking the order of the stabilizers given in Lemma~\ref{Stabilizer},
we have $N\cong\ZZ_2^r$, $\ZZ_3^l$ or $\ZZ_5^k$, where $1\le r\le 10$, $1\le
l\le3$ and $1\le k\le2$. By Lemma~\ref{N/C},
$F/C_F(N)\lesssim \Aut(N)\cong \GL(r,2)$, $\GL(l,3)$ or $\GL(k,5)$.
Clearly, $N\le C_F(N)$. If $N=C_F(N)$, then $\A_{21}\cong S\cong F/N=F/C_F(N)
\lesssim \GL(r,2)$, $\GL(l,3)$ or $\GL(k,5)$. However, by
Magma \cite{Magma}, each of $\GL(r,2)$, $\GL(l,3)$ and $\GL(k,5)$ has no subgroup
isomorphic to $\A_{21}$ for $1\le r\le 10$, $1\le
l\le3$ and $1\le k\le2$, a contradiction.
Hence $N<C_F(N)$ and $1\not=C_F(N)/N\lhd F/N\cong\A_{21}$. It
follows $F=C_F(N)=N\times S$,
$F_\a/S_\a\cong F/S\cong N$,
and $F_\a$ is soluble because $S_{\a}\cong\F_{21}$.
By Lemma~\ref{Stabilizer}, we conclude
$F_\a\cong \F_{42}$, $\F_{42}{\times}\ZZ_2$ or $\F_{21}\times\ZZ_3$.
A direct computation by Magma \cite{Magma} shows
there is no feasible element to $F$ and $F_{\a}$,
it is also a contradiction.

Thus, $\A$ is quasiprimitive on $V\Ga$.
Let $M$ be a minimal normal subgroup of $\A$. Then $M=D^d$,
with $D$ a non-abelian simple group, is transitive on
$V\Ga$,
so $|V\Ga|=|\A_{20}|$ divides $|M|$ and $19\div|D|$.
If $d\ge2$, then $19^2\div|M|$,
which is a contradiction because $|\A|\div|\A_{20}|\cdot2^{24}\cdot3^4\cdot5^2\cdot7$
is not divisible by $19^2$.
Hence $d=1$ and $M=D\lhd \A$. Let $C=C_\A(D)$. Then $C\lhd \A$
and $CD = C{\times}D$. If $C\not=1$, then $C$ is transitive on $V\Ga$ as
$\A$ is quasiprimitive on $V\Ga$,
with similar discussion as above, we have $C$ is insoluble and $19\div|C|$.
Therefore, $19^2\div|CD|$, again a contradiction.
Hence $C=1$ and $\A$ is almost simple.

Since $M\cap S\unlhd S\cong\A_{21}$,
$M\cap S=1$ or $S$. If $M\cap S=1$, then
$|M|\div 2^{24}\cdot 3^3\cdot5^2$, it is a contradiction as $|\A_{20}|\div|M|$.
Thus, $M\cap S=S$ and so $S\le M$. Hence $M$ is a non-abelian simple
group satisfying $|\A_{21}|\div|M|\div|\A_{20}|\cdot2^{24}\cdot3^4\cdot5^2\cdot7$.
By \cite[P.135-136]{Gorenstein}, we conclude $M=S\cong\A_{21}$.
Thus $\A\le\Aut(M)\cong\S_{21}$. If $\A\cong\S_{21}$, then $|\A_\a|=|\A:T|=42$,
and so  $\A_\a\cong\F_{42}$ by Lemma \ref{Stabilizer}. A direct computation by Magma \cite{Magma} shows
there is no feasible element to $\A$ and $\A_{\a}$,
a contradiction. Hence $\A\cong\A_{21}$
and $\Ga$ is $1$-transitive.
\qed


\begin{thebibliography}{99}

\bibitem{Magma}
W. Bosma, J. Cannon, C. Playoust, The MAGMA algebra system I:
The user language, {\it J. Symbolic Comput.} \textbf{24} (1997), 235--265.

\bibitem{Atlas}
J. H. Conway, R. T. Curtis, S. P. Norton, R. A. Parker,
R. A. Wilson, {\it Atlas of Finite Groups},
Oxford Univ. Press, London/New York, 1985.

\bibitem{DM}
J. D. Dixon, B. Mortimer, {\it Permutation groups}, Sprimger-Verlag, 1996.

\bibitem{DFZ17}
J. L. Du, Y. Q. Feng, J. X. Zhou,
Pentavalent symmetric graphs admitting vertex-transitive non-abelian simple groups,
{\it Europ. J. Combin.} \textbf{63} (2017), 134--145.



\bibitem{FLX04}
X. G. Fang, C. H. Li, M. Y. Xu,
On edge-transitive Cayley graphs of valency four, {\it Europ. J. Combin.}
{\bf 25} (2004), 1107--1116.

\bibitem{FMW11}
X. G. Fang, X. S. Ma, J. Wang, On locally-primitive Cayley graphs of finite simple groups,
{\it J. Combin. Theory Ser. A}
{\bf 118} (2011), 1039--1051.

\bibitem{FPW02}
X. G. Fang, C. E. Praeger, J. Wang, On the automorphism group of Cayley graphs of finite simple groups,
{\it J. London Math. Soc.}
{\bf 66} (2002), 563--578.



\bibitem{Gorenstein}
D. Gorenstein, {\it Finite Simple Groups}, Plenum Press, New
York, 1982.


\bibitem{GLH16}
S. T. Guo, Y. T. Li, X. H. Hua, $(G,s)$-Transitive Graphs of Valency 7,
\emph{Algebra Coll.} {\bf 23} (2016), 493--500.

\bibitem{Huppert}
B. Huppert, {\it Finite Groups}, Springer-Verlag, Berlin, 1967.


\bibitem{Huppert-L}
B. Huppert, W. Lempken, Simple groups of order divisible by at most four primes,
{\it Proc. F. Scorina Gomel State Univ.}
{\bf 10} (2000), 64--75.

\bibitem{KL90}
P. Kleidman, M. Liebeck, {\it The Subgroup Structure of The Finite Classical Groups},
Cambridge Univ. Press, 1990.

\bibitem{Li96}
C. H. Li, Isomorphisms of finite Cayley graphs,
Ph. D. Thesis,
The University of Western Australia, 1996.

\bibitem{Li-ETCay}
C. H. Li, Finite edge-transitive Cayley graphs and rotary Cayley
maps, {\it Trans. Amer. Math. Soc.} {\bf 358} (2006), 4605--4635.



\bibitem{LLW16}
C. H. Li, Z. P. Lu, G. X. Wang,
Arc-transitive graphs of square-free order and small valency,
\emph{Discrete Math.} {\bf 339} (2016), 2907--2918.

\bibitem{Li}
C. H. Li, J. M. Pan, Finite 2-arc-transitive abelian Cayley graphs,
{\it Europ. J. Combin.} \textbf{29} (2008), 148--158.

\bibitem{MR90}
B. D. McKay, G. Royle, The transitive graphs with at most 26 vertices,
{\it Ars Combin.} {\bf 30} (1990), 161--176.

\bibitem{Praeger92}
C. E. Praeger. An O'Nan-Scott theorem for finite quasiprimitive
permutation groups and an application to 2-arc transitive graphs.
{\it J. London. Math. Soc.} {\bf 47} (1993), 227--239.

\bibitem{Sab64}
B. O. Sabidussi, Vertex-transitive graphs,
{\it Monash Math.} {\bf 68} (1964), 426--438.

\bibitem{Schur}
I. Schur, Untersuchen \" uber die Darstellung der endlichen Gruppen durch gebrochenen linearen Substitutionen.
{\it J. Reine Angew. Math.} {\bf 132} (1907), 85--137.

\bibitem{Suzuki}
M. Suzuki, {\it Group Theroy II}, Springer-Verlag, New York, 1985.

\bibitem{Xu05}
S. J. Xu, X. G. Fang, J. Wang, M. Y. Xu,
On cubic s-arc-transitive Cayley graphs on finite simple groups,
{\it Europ. J. Combin.} {\bf 26} (2005), 133--143.

\bibitem{Xu07}
S. J. Xu, X. G. Fang, J. Wang, M. Y. Xu,
5-arc-transitive cubic Cayley graphs on finite simple groups,
{\it Europ. J. Combin.} {\bf 28} (2007), 1023--1036.


\end{thebibliography}
\end{document}